\documentclass[11pt,fleqn,twoside]{article}
\usepackage{amsfonts,latexsym}
\makeatletter
\newcommand{\prava}{\footnotesize\it
\begin{flushright}
\begin{minipage}{18cm}
Copyright \copyright 1998 by M. Wang, Y. Zhou and Z. Li
\end{minipage}
\end{flushright}}

\newcommand{\name}[1]{\begin{flushleft}
                       \LARGE \bf #1
                       \end{flushleft}\vspace{-3mm}}

\newcommand{\Author}[1]{\begin{flushleft}
                       \it #1 \end{flushleft}}

\newcommand{\Adress}[1]{\begin{flushleft}
                       \it #1 \end{flushleft}}

\newcommand{\Date}[1]{\begin{flushleft}
                      \small  \it #1 \end{flushleft}}

\newcommand{\ehkol}{Author \ name}
\newcommand{\ohkol}{Article \ name}
\renewcommand{\@evenhead}{
\hspace*{-3pt}\raisebox{-15pt}[\headheight][0pt]{\vbox{\hbox to \textwidth
{\thepage \hfil \ehkol}\vskip4pt \hrule}}}
\renewcommand{\@oddhead}{
\hspace*{-3pt}\raisebox{-15pt}[\headheight][0pt]{\vbox{\hbox to \textwidth
{\ohkol \hfil \thepage}\vskip4pt\hrule}}}
\renewcommand{\@evenfoot}{}
\renewcommand{\@oddfoot}{}

\newcommand{\be}{\begin{equation}}
\newcommand{\ee}{\end{equation}}
\newcommand{\ba}{\hspace*{-5pt}\begin{array}}
\newcommand{\ea}{\end{array}}

\newcommand{\ds}{\displaystyle}

\makeatother

\begin{document}
\setcounter{page}{120}

\thispagestyle{empty}

\renewcommand{\ehkol}{M. Wang, Y. Zhou and Zh. Li}
\renewcommand{\ohkol}{A Nonlinear
Transformation of the Dispersive Long Wave Equations}

\begin{flushleft}
\footnotesize \sf
Journal of Nonlinear Mathematical Physics \qquad 1998, V.5, N~2,
\pageref{wang-fp}--\pageref{wang-lp}. \hfill {\sc Letter}
\end{flushleft}

\vspace{-5mm}

\renewcommand{\footnoterule}{}
{\renewcommand{\thefootnote}{}
 \footnote{\prava}}

\name{A Nonlinear Transformation of the Dispersive Long Wave Equations
in (2+1) Dimensions \\ and its Applications}\label{wang-fp}

\Author{Mingliang WANG, Yubin ZHOU and Zhibin LI}

\Adress{Department of Mathematics, Lanzhou University,\\
  Lanzhou, Gansu, 730000, P.R. of China}

\Date{Received January 23, 1998; Accepted March 12, 1998}

\begin{abstract}
\noindent
 A nonlinear transformation of the dispersive long wave
equations in ($2+1$) dimensions is derived by using the homogeneous
balance
method. With the aid of the transformation given here, exact solutions
of the equations are obtained.
\end{abstract}

\section{Introduction: New results}\label{wang:sec1}

\setcounter{equation}{0}
\renewcommand{\theequation}{\arabic{section}.\arabic{equation}}

 Consider the dispersive long wave equations in ($2+1$) dimensions,
\be\label{wang:eq1.1}
u_{yt} + h_{xx} + {1 \over 2} \left(u^2\right)_{xy} = 0,
\ee
\be  \label{wang:eq1.2}
 h_t + (u h + u + u_{xy})_x = 0
\ee
which were f\/irst obtained by Boiti {\it et al}~\cite{wang:ref1} as a
compatibility condition for a ``weak" Lax pair. A Kac-Moody-Virasoro
type Lie algebra for eqs. (\ref{wang:eq1.1}) and (\ref{wang:eq1.2})
were given by Paquin and Winternitz~\cite{wang:ref2}. Moreover, Sen-yue
Lou~\cite{wang:ref3} showed that eqs. (\ref{wang:eq1.1}) and
(\ref{wang:eq1.2}) do not pass the Painlev\'e test, both in the ARS
algorithm and in the WTC approach.  Equations (\ref{wang:eq1.1}) and
(\ref{wang:eq1.2}) can be reduced to the ($1+1$) dimensional
model~\cite{wang:ref4}
\be\label{wang:eq1.3}
 u_t + h_z + {1 \over 2} \left(u^2\right)_z = 0,
 \ee
\be  \label{wang:eq1.4}
 h_t + (u h + u + u_{zz})_z = 0,
\ee
for $u = u(x+y, t) = u(z, t)$, $h = h(x+y, t) = h(z, t)$; the solitary
wave solutions of which were
obtained by the f\/irst author~\cite{wang:ref5}.

In the present paper we  show that if $\varphi = \varphi(x,y,t)$ is a
nontrivial solution of the linear equation
\begin{equation}\label{wang:eq1.5}
 \varphi_t \pm \varphi_{xx} = 0,
\end{equation}
then the nonlinear transformation
\begin{equation}\label{wang:eq1.6}
 u(x,y,t) = \pm \frac{2\varphi_x}{\varphi},
 \qquad
 h(x,y,t) = - \frac{2\varphi_x\varphi_y}{\varphi^2}
 + \frac{2\varphi_{xy}}{\varphi} -1,
\end{equation}
satisf\/ies eqs. (\ref{wang:eq1.1}) and (\ref{wang:eq1.2}).

 At the time of writing, we have not seen the results like these in the
literature. If this assertion is true, then the exact
solutions of eqs. (\ref{wang:eq1.1}) and (\ref{wang:eq1.2}) can be
easily obtained by the use of this transformation.
In fact, we know that linear equation
(\ref{wang:eq1.5}) can be handled by a variety of methods, and has many
nontrivial solutions which may include arbitrary
functions of $y$, so do eqs. (\ref{wang:eq1.1}) and
(\ref{wang:eq1.2}) according to transformation (\ref{wang:eq1.6}).
For example, the function
\begin{equation}\label{wang:eq1.7}
 \varphi = 1 + \exp\left[a(y) x \mp a^2(y) t + b(y)\right],
\end{equation}
satisf\/ies eq. (\ref{wang:eq1.5}), where $a=a(y)$ and $b=b(y)$ are
both arbitrary dif\/ferentiable functions. Thus substituting
(\ref{wang:eq1.7}) into (\ref{wang:eq1.6}), we obtain, after doing some
computations, the exact solutions of
eqs. (\ref{wang:eq1.1}) and (\ref{wang:eq1.2}) as follows:
\be\label{wang:eq1.8}
\ds   u(x,y,t)  =  \pm a(y)
\left\{1+ \tanh\frac{1}{2}\left[a(y) x \mp a^2(y) t +
b(y)\right]\right\},
\ee
\be  \label{wang:eq1.9}
\ba{l}
\ds   h(x,y,t)  =  \frac{1}{2} a(y)
\left[a^\prime(y) x \mp 2a(y) a^\prime(y) t +
b^\prime(y)\right]
    \mbox{sech}^2 \frac{1}{2}
    \left[ a(y) x \mp a^2(y)t + b(y)\right]  \\[2mm]
\ds \qquad \qquad + a^\prime(y)\tanh\frac{1}{2}
\left[a(y) x \mp a^2(y) t + b(y)\right] + a^\prime(y) -1,
\ea\!\!
\ee
where it is assumed that $\lim\limits_{y\rightarrow \pm \infty} a(y) =
\mbox{const.}$; $\lim\limits_{y \rightarrow \pm \infty} b^\prime(y)
=\mbox{const.}$ so that $u$ and $h$ are bounded as $y$ goes to
inf\/inity.

In particular, if $a = \mbox{const}$, $b = c y + d$, $c$ and $d$ are
both arbitrary constants, then (\ref{wang:eq1.8}) and
(\ref{wang:eq1.9}) become

\setcounter{equation}{7}
\renewcommand{\theequation}{\arabic{section}.\arabic{equation}$'$}
\be\label{wang:eq1.8p}
 u(x,y,t) = \pm a \left[1 + \tanh\frac{1}{2}
 \left(ax \mp a^2 t + c y + d\right)\right],
 \ee
 \be  \label{wang:eq1.9p}
h(x,y,t) = \frac{ac}{2} \mbox{sech}^2 \frac{1}{2}
\left[a x \mp a^2 t + c y + d\right]- 1,
\ee
which coincide with the results in Ref.~\cite{wang:ref6} (however, only one
solution was given there). If taking $a = c$ and $x+y = z$, then
(\ref{wang:eq1.8p}) and (\ref{wang:eq1.9p}) become the results in
Ref.\cite{wang:ref5} (also, only one solution was given there).

\setcounter{equation}{0}
\renewcommand{\theequation}{\arabic{section}.\arabic{equation}}

\section{Derivation of the nonlinear transformation}\label{wang:sec2}

Now we use the homogeneous balance method to derive the nonlinear
transformation (\ref{wang:eq1.6}). A function $\varphi =
\varphi(x,y,t)$ is called a quasisolution of eqs. (\ref{wang:eq1.1})
and (\ref{wang:eq1.2}), if there exists a pair of functions $f =
f(\varphi)$ and $g=g(\varphi)$ of one variable only, so that the
following expressions
\be\label{wang:eq2.1}
\ba{l}
\ds u(x,y,t) = \frac{ \partial^{l + m +n} f(\varphi) }{
              \partial x^l \partial y^m \partial t^n }
    + \mbox{a suitable linear combination of some partial}\\[3mm]
\qquad\qquad\qquad    \mbox{derivatives  (their
orders are all lower than} \ l+m+n)\\[1mm]
\qquad\qquad\qquad \mbox{with respect to } x, y \mbox{ and } t \mbox{
of } f(\varphi),
\ea
\ee
\be  \label{wang:eq2.2}
 h(x,y,t) = \frac{\partial^{p+q+r} g(\varphi)}{
      \partial x^p \partial y^q \partial t^r} + \cdots,
\ee
which can be rewritten as
\setcounter{equation}{0}
\renewcommand{\theequation}{\arabic{section}.\arabic{equation}$'$}
\be\label{wang:eq2.1p}
\ba{l}
 u(x,y,t) = f^{(l + m+ n)}(\varphi) \varphi_x^l \varphi_y^m \varphi_t^n +
  \mbox{a polynomial of various partial}\\[1mm]
\qquad \qquad \qquad \mbox{derivatives of }  \varphi(x,y,t)
(\mbox{ in spite of } f(\varphi) \ \mbox{and its}\\[1mm]
\qquad \qquad \qquad \mbox{derivatives,
   the degree of which} \mbox{ is lower than } l+m+n),
   \ea
 \ee
 \be \label{wang:eq2.2p}
 h(x,y,t) = g^{(p+q+r)}(\varphi) \varphi_x^p
 \varphi_y^q \varphi_t^r + \cdots
\ee
is actually a solution of eqs. (\ref{wang:eq1.1}) and
(\ref{wang:eq1.2}). Here $f = f(\varphi)$, $g = g(\varphi)$,
$\varphi = \varphi(x,y,t)$, and nonnegative integers $l, m, n$
and $p, q, r$ are all to be determined later.  The unwritten
parts in (\ref{wang:eq2.2}) and (\ref{wang:eq2.2p}) are similar to
those in (\ref{wang:eq2.1}) and (\ref{wang:eq2.1p}), respectively. For
the sake of simplicity, we
omit these similar parts without explanation.

\setcounter{equation}{2}
\renewcommand{\theequation}{\arabic{section}.\arabic{equation}}

 The homogeneous balance method mainly consists of three steps:

{\bf  First step:} Determine nonnegative integers $l, m, n$ and $p, q, r$
which determine the order of the highest (partial) derivative of
$f(\varphi)$ and $g(\varphi)$ in (\ref{wang:eq2.1}) and
(\ref{wang:eq2.2}), respectively. Then
\be\label{wang:eq2.3}
(u^2)_{xy} =
\left\{\left[f^{(l+m+n)}(\varphi)\right]^2\right\}^{\prime\prime}
    \varphi_x^{2l+1} \varphi_y^{2m+1} \varphi_t^{2n} + \cdots,
\ee
\be \label{wang:eq2.4}
 h_{xx} = g^{(p+q+r+2)}(\varphi) \varphi_x^{p+2}
 \varphi_y^{q} \varphi_t^{r} +\cdots,
\ee
\be \label{wang:eq2.5}
(u h)_x = \left[f^{(l+m+n)}(\varphi) g^{(p+q+r)}(\varphi)\right]^\prime
   \varphi_x^{l+p+1} \varphi_y^{m+q} \varphi_t^{n+r} + \cdots,
\ee
\be \label{wang:eq2.6}
u_{xxy} = f^{(l+m+n+3)}(\varphi) \varphi_x^{l+2}
\varphi_y^{m+1} \varphi_t^{n} + \cdots.
\ee
In order that $(u^2)_{xy}$ and $h_{xx}$ in eq. (\ref{wang:eq1.1}) can be
partially balanced in view of (\ref{wang:eq2.3}) and
(\ref{wang:eq2.4}), it should be required that
the highest degrees in partial derivative for $\varphi(x,y,t)$ appearing in
$(u^2)_{xy}$ and $h_{xx}$ are equal; and similarly, in order that $(uh)_x$
and $u_{xxy}$ in eq.~(\ref{wang:eq1.2}) can be partially balanced in view of
(\ref{wang:eq2.5}) and (\ref{wang:eq2.6}). It should also be required
that the highest degrees in partial derivatives for $\varphi(x,y,t)$
appearing in $(u h)_x$ and $u_{xxy}$ are equal. Thus we have
\begin{equation}\label{wang:eq2.7}
\ba{l}
 \left.
  \ba{l}
    2l + 1 = p+2,\ 2m+1 = q,\ 2n = r \\
    l+p+1 = l+2,\  m+q = m+1,\ n+r = n
  \ea
 \right\}\\[4mm]
\qquad  \qquad \qquad \qquad \qquad
  \Longrightarrow l =1,\ m=n=0;\ p=q=1,\ r=0 .
\ea
\end{equation}
If the results in (\ref{wang:eq2.7}) are used, then the expressions
(\ref{wang:eq2.1}) and (\ref{wang:eq2.2}) (or (\ref{wang:eq2.1p}) and
(\ref{wang:eq2.2p})) are of the form
\begin{equation}\label{wang:eq2.8}
 u = f_x(\varphi) = f^\prime \varphi_x, \qquad
 h = g_{xy}(\varphi) + A = g^{\prime\prime} \varphi_x \varphi_y +
 g^\prime \varphi_{xy} + A,
\end{equation}
where $A$ is a constant to be determined.

{\bf Second step} Determine functions $f = f(\varphi)$ and $g = g(\varphi)$
in expressions (\ref{wang:eq2.1}) and (\ref{wang:eq2.2})
((\ref{wang:eq2.1p}) and (\ref{wang:eq2.2p})).
From (\ref{wang:eq2.8}) it is easy to deduce that
\be\label{wang:eq2.9}
 u_{yt} = f^{\prime\prime\prime} \varphi_x \varphi_y \varphi_t
   + f^{\prime\prime} (\varphi_{xt} \varphi_y + \varphi_x \varphi_{yt} +
   \varphi_{xy} \varphi_t)    + f^\prime \varphi_{xyt},
\ee
\be \label{wang:eq2.10}
\ba{l}
h_{xx} = g^{(4)} \varphi_x^3 \varphi_y
   + g^{\prime\prime\prime}\left( 3 \varphi_x \varphi_{xx} \varphi_y +
   3 \varphi_x^2 \varphi_{xy}\right)\\[2mm]
\qquad \qquad
   + g^{\prime\prime}(\varphi_{xxy} \varphi_x +
   3 \varphi_{xx} \varphi_{xy} + 3 \varphi_x\varphi_{xxy})
   + g^\prime \varphi_{xxxy},
\ea
\ee
\be \label{wang:eq2.11}
\ba{l}
\frac{1}{2}(u^2)_{xy} = (f^{\prime\prime 2} + f^\prime
f^{\prime\prime\prime}) \varphi_x^3 \varphi_y \\[2mm]
\qquad \qquad+ f^\prime f^{\prime\prime}(3 \varphi_x^2
\varphi_{xy} + 2 \varphi_x \varphi_{xx} \varphi_y)
+ f^{\prime 2} (\varphi_{xx} \varphi_{xy} + \varphi_x \varphi_{xxy}),
\ea
\ee
\be \label{wang:eq2.12}
 h_t = g^{\prime\prime\prime} \varphi_x \varphi_y \varphi_t
   + g^{\prime\prime} (\varphi_{xt} \varphi_y + \varphi_{yt} \varphi_x +
   \varphi_{xy} \varphi_t)    + g^\prime \varphi_{xyt},
\ee
\be\label{wang:eq2.13}
\ba{l}
(u h)_x = (f^{\prime\prime} g^{\prime\prime} +
f^\prime g^{\prime\prime\prime}) \varphi_x^3 \varphi_y
+ f^\prime g^{\prime\prime} (2 \varphi_x \varphi_{xx} \varphi_y +
\varphi_x^2 \varphi_{xy})
+ f^{\prime\prime} g^\prime \varphi_x^2\varphi_{xy}
 \\[2mm]
\qquad \qquad+ f^\prime g^{\prime\prime} \varphi_x^2
\varphi_{xy}  + f^\prime g^\prime (\varphi_{xx} \varphi_{xy}
+ \varphi_x \varphi_{xxy})
   + A f^{\prime\prime} \varphi_x^2 + A f^\prime \varphi_{xx},
\ea
\ee
\be \label{wang:eq2.14}
 u_x = f^{\prime\prime} \varphi_x^2 + f^\prime \varphi_{xx},
\ee
\be \label{wang:eq2.15}
\ba{l}
(u_{xy})_x = f^{(4)} \varphi_x^3 \varphi_y  + f^{\prime\prime\prime}
( 3 \varphi_x^2 \varphi_{xy} + 3 \varphi_x \varphi_{xx} \varphi_y)
  \\[2mm]
  \qquad \qquad+ f^{\prime\prime}( 3 \varphi_{xx} \varphi_{xy} +
  3 \varphi_x \varphi_{xxy} + \varphi_{xxx}\varphi_y)
  + f^\prime \varphi_{xxxy} .
  \ea
\ee
Substituting (\ref{wang:eq2.9})--(\ref{wang:eq2.11}) and
(\ref{wang:eq2.12})--(\ref{wang:eq2.15}) into the left hand sides of
eqs. (\ref{wang:eq1.1}) and (\ref{wang:eq1.2}), respectively, and
collecting all homogeneous terms of 4, 3, 2 and 1 degree in partial
derivatives of $\varphi(x,y,t)$ together, yields
\be\label{wang:eq2.16}
\ba{l}
 u_{yt} + h_{xx} + \frac{1}{2} (u^2)_{xy} =
   (g^{(4)} + f^{\prime\prime 2} +
   f^\prime f^{\prime\prime\prime}) \varphi_x^3 \varphi_y \\[2mm]
\qquad + [f^{\prime\prime\prime} \varphi_x \varphi_y \varphi_t
     + g^{\prime\prime\prime} (3 \varphi_x \varphi_{xx} \varphi_y +
     3 \varphi_x^2 \varphi_{xy})  + f^\prime f^{\prime\prime}
      (3 \varphi_x^2 \varphi_{xy} + 2 \varphi_x \varphi_{xxy})]
 \\[2mm]
\qquad  + [f^{\prime\prime} (\varphi_{xt} \varphi_y +
\varphi_{yt} \varphi_x + \varphi_{xy} \varphi_t)
+ g^{\prime\prime} (\varphi_{xxx} \varphi_y +
3 \varphi_x \varphi_{xxy} + 3 \varphi_{xx}\varphi_{xy})
\\[2mm]
\qquad + f^{\prime 2} (\varphi_{xx} \varphi_{xy} + \varphi_x \varphi_{xxy})]
   + (f^\prime \varphi_{xyt} + g^\prime \varphi_{xxxy} )
\ea
\ee
and
\be\label{wang:eq2.17}
\ba{l}
h_t + (u h + u + u_{xy})_x = ( f^{(4)} + f^{\prime\prime} g^{\prime\prime}
+ f^\prime g^{\prime\prime\prime}) \varphi_x^3 \varphi_y
\\[2mm]
\qquad  + [ g^{\prime\prime\prime} \varphi_x \varphi_y \varphi_t
   + f^\prime g^{\prime\prime} (2 \varphi_x \varphi_{xx} \varphi_y
   + \varphi_x^2 \varphi_{xy}) + ( f^{\prime\prime} g^\prime
   + f^\prime g^{\prime\prime}) \varphi_x^2 \varphi_{xy}
\\[2mm]
\qquad   + f^{\prime\prime\prime} ( 3 \varphi_x^2 \varphi_{xy} +
3 \varphi_x \varphi_{xx} \varphi_y)]
  + [g^{\prime\prime} (\varphi_{xt} \varphi_y +
  \varphi_{yt} \varphi_x + \varphi_{xy} \varphi_t)
 \\[2mm]
\qquad + g^\prime f^\prime ( \varphi_{xx} \varphi_{xy} +
\varphi_x \varphi_{xxy})  + f^{\prime\prime}
(3 \varphi_{xx} \varphi_{xy} + 3 \varphi_x \varphi_{xxy}
    + \varphi_{xxx} \varphi_y) \\[2mm]
\qquad  + (A+1) f^{\prime\prime} \varphi_x^2 ]
+ [ g^\prime \varphi_{xyt} +
(A + 1) f^\prime \varphi_{xx} + f^\prime \varphi_{xxxy}].
\ea
\ee
Setting the coef\/f\/icients of $\varphi_x^3 \varphi_y$ in
(\ref{wang:eq2.16}) and (\ref{wang:eq2.17}) to zero, yields a system of
ordinary dif\/ferential equations for $f(\varphi)$ and
$g(\varphi)$, namely
\begin{equation}\label{wang:eq2.18}
 \left.
 \ba{l}
  g^{(4)} + f^{\prime\prime 2} + f^\prime f^{\prime\prime\prime} = 0,\\[1mm]
  f^{(4)} + f^{\prime\prime} g^{\prime\prime} +
  f^\prime g^{\prime\prime\prime} = 0,
 \end{array}
 \right\}
\end{equation}
which admits two solutions
\begin{equation}\label{wang:eq2.19}
 f(\varphi) = \pm 2 \ln \varphi, \qquad g = 2\ln \varphi ,
\end{equation}
and thereby
\begin{equation}\label{wang:eq2.20}
 g^\prime g^{\prime\prime} = - g^{\prime\prime\prime},
 \qquad g^{\prime 2} = - 2 g^{\prime\prime} .
\end{equation}

{\bf  Third step:} Determine the equation satisf\/ied by the quasisolution
$\varphi = \varphi (x, y, t)$ of eqs. (\ref{wang:eq1.1})
and (\ref{wang:eq1.2}).
Using (\ref{wang:eq2.18})--(\ref{wang:eq2.20}), the expressions (\ref{wang:eq2.16}) and
(\ref{wang:eq2.17}) can be simplif\/ied as
\be\label{wang:eq2.21}
\ba{l}
\ds  u_{yt} + \eta_{xx} + \frac{1}{2} (u^2)_{xy}\\[2mm]
\ds   \qquad \qquad =
   \pm \left[\varphi_x \varphi_y g^{\prime\prime\prime}
   + g^{\prime\prime}\left(\varphi_x \frac{\partial}{\partial y}
   + \varphi_y \frac{\partial}{\partial x} + \varphi_{xy}\right)
   + f^\prime \frac{\partial^2}{\partial x \partial y}\right] (\varphi_t \pm
\varphi_{xx}),
\ea
\ee
\be\label{wang:eq2.22}
\ba{l}
\ds \eta_t + (u \eta + u + u_{xy})_x \\[2mm]
\ds  \qquad \qquad =
\left[ \varphi_x \varphi_y g^{\prime\prime\prime}
   + g^{\prime\prime}\left(\varphi_x \frac{\partial}{\partial y}
   + \varphi_y \frac{\partial}{\partial x} + \varphi_{xy}\right)
   + g^\prime \frac{\partial^2}{\partial x \partial y}\right]
   (\varphi_t \pm \varphi_{xx}),
   \ea
\ee
provided that we take
\begin{equation}\label{wang:eq2.23}
 A = - 1 .
\end{equation}

 In view of (\ref{wang:eq2.21}) and (\ref{wang:eq2.22}),
it is easily deduced that if
$\varphi = \varphi(x,y,t)$ satisf\/ies the linear equation
\[
 \varphi_t \pm \varphi_{xx} = 0 ,
\]
i.e., eq. (\ref{wang:eq1.5}), then the right
hand sides of (\ref{wang:eq2.21}) and
(\ref{wang:eq2.22}) vanish respectively. This means
that (\ref{wang:eq2.8}) actually
solves eqs. (\ref{wang:eq1.1}) and (\ref{wang:eq1.2}), provided that
\[
 f = \pm 2 \ln \varphi, \qquad g = 2 \ln \varphi ,
 \qquad (\mbox{ i.e., } (\ref{wang:eq2.19}))
\]
and
\[
 A = - 1 . \qquad (\mbox{ i.e., } (\ref{wang:eq2.23}))
\]
Substituting (\ref{wang:eq2.19}) and (\ref{wang:eq2.23}) into (\ref{wang:eq2.8}), yields
\[
 u = \pm \frac{2 \varphi_x}{\varphi},
 \qquad  \eta = - \frac{2 \varphi_x \varphi_y}{\varphi^2}
 + \frac{2 \varphi_{xy}}{\varphi} - 1 ,
\]
which is the nonlinear transformation (\ref{wang:eq1.6}).

\subsection*{Acknowledgment}

This work was supported by the National Natural Science Foundation of
China.

\label{wang-lp}
\end{document}